\documentstyle[12pt]{article}
\newcommand{\R}{\mbox{$I \kern -4pt R$}}             
\oddsidemargin = - 0.5cm
\begin {document}
\baselineskip 15pt
\title{Thermal convection in a cylindrical annulus heated laterally} 
 \author{S. Hoyas$\star$, H. Herrero$\star$ and A. M. Mancho$_\dagger$\\
     {$\star$ \small Departamento de Matem\'aticas}, 
 	{\small Facultad de Ciencias Qu\'{\i}micas}\\
        {\small Universidad de Castilla-La Mancha}, 
        {\small 13071 Ciudad Real, Spain.}\\
        {$\dagger$ \small Department of Mathematics}, 
{\small School of Mathematics}\\
     {\small University of Bristol},    
{ \small University Walk,} \\
{ \small  Bristol BS8 1TW, United Kingdom.}\\
{\small E-mail: {\tt Sergio.Hoyas}@{\tt uclm.es} (S Hoyas), {\tt  
Henar.Herrero@uclm.es}}\\
{\small  (H Herrero) and {\tt A.M.Mancho@Bristol.ac.uk} (A M  
Mancho) }\\}
\maketitle
\noindent

\section*{Abstract}

In this paper we study
thermoconvective instabilities appearing in a fluid within a 
cylindrical annulus heated laterally.  
As soon as a horizontal temperature gradient is applied a 
convective state appears. As the temperature gradient reaches
a critical value a stationary or oscillatory bifurcation may take place.
The problem is modelled with a novel method which extends the one
described in \cite{numerico}. The Navier Stokes equations are
solved in the primitive variable  formulation,  with appropriate boundary  
conditions for pressure. This is a low order formulation which in
cylindrical coordinates introduces lower order singularities. The
problem is discretized with a Chebyshev  
collocation method easily implemented and its convergence
has been checked. The results obtained are not only  in very
good agreement with those obtained in experiments, but also provide a deeper
insight into important physical parameters developing the instability,
which has not been reported before.
 
\section{Introduction}

The problem of thermoconvective instabilities in fluid layers driven 
by a temperature gradient 
has become a classical subject in fluid mechanics \cite{benard,pearson}.   
It is well known
that two different effects are responsible for the onset of motion when   
the
temperature difference becomes larger than a certain threshold: 
gravity and capillary forces. When both effects are 
taken into account the problem  
is called 
B\'enard-Marangoni (BM) convection \cite{nield}.\par 
To solve numerically hydrodynamical problems  in the
primitive variable formulation raises two questions, one is the 
handling of first order derivatives for  pressure 
and the other is finding its boundary conditions  
\cite{gresho,poz,orszag}, which are not intuitive. 
In thermoconvective problems pressure usually is avoided, 
for instance in Refs. \cite{mer1,mer2}
 the method of potentials of velocity is used
to eliminate it from the equations. This 
technique raises the order of
the differential equations and additional boundary conditions may be
 required. This is particularly troublesome in cylindrical coordinates
 where high order derivatives cause awkward difficulties. In Ref. 
\cite{dauby2} the primitive variables formulation is used
although the use of a particular spectral method allows the removal of  pressure. 
In Ref. \cite{numerico} the linear stability analysis of two convection
 problems are solved in the primitive variables formulation  and taking  appropriate  
boundary conditions for pressure. In this article very good accuracy is  obtained
both for cartesian and cylindrical coordinates. This method has been  succesfully
applied to describe experimental results in cylindrical containers  \cite{mluisa}. 
In this paper 
we extend this   method  to study a BM problem in a cylinder heated  laterally. 

The physical set up (Fig. 1) consists of a fluid filling up a 
container bounded by two concentric cylinders. 
The upper surface is open  
to the air and the fluid is heated laterally through the lateral 
walls which are conducting.
This set up corresponds to the experiment 
reported in Ref. \cite{angel}, and also it is 
comparable to similar experiments in rectangular containers \cite{daviaud,burguete}. In \cite{arancha}
a similar problem is treated  theoretically, however  in this work  
rotation is considered and the gravity field is along
the radial coordinate. In the  problem we  treat,
as soon as a slight
difference of temperature is imposed between the lateral walls, a  
stationary solution appears which is called
basic state. When the temperature of the walls 
is modified the basic state can become unstable and  
bifurcate to different patterns, both  stationary and oscillatory.
When the bifurcation is stationary, the  emerging 
pattern consists of a 3D structure of rolls whose axes are parallel 
to the radial coordinate (radial rolls) and in the oscillatory case,
the axes are tilted with respect the radial axis. These results 
coincide with those obtained experimentally \cite{angel,daviaud,burguete}.

The basic state is calculated solving a set of nonlinear equations  whereby 
the first contribution is obtained by a linear approach. To improve 
the solution 
we expand the corrections of the unknown fields in Chebyshev  polynomials, 
and we pose the equations at the Gauss-Lobatto 
collocation points  \cite{canuto}. 
The bifurcation thresholds are obtained through a generalized eigenvalue  
problem with the same Chebyshev collocation method. 
The convergence of the method is studied by
comparing different expansions.

The organization of the article is
as follows.  In Sec. 2 the 
formulation of the problem with the equations and boundary  
conditions is explained. In the third section the
basic state is calculated and the results obtained are discussed
for different physical conditions. In the fourth  
section the linear stability of the basic solution 
is performed and instabilities are studied for different
parameters. In the fifth section conclusions are detailed.

\section{Formulation of the problem}
The physical set up considered is shown in Fig. 1. 
A horizontal fluid  layer  of depth $d$ ($z$ coordinate) is in 
a  container bounded by  
two concentric cylinders of radii $a$ and $a+\delta $ ($r$ coordinate).  The
bottom plate is rigid and the top is open to the
atmosphere. The inner cylinder has a temperature $T_{\max }$  
whereas the outer one is at $T_{\min}$ and the environment 
is at $T_0$. We define
$\triangle T = T_{\max}-T_{0}$ and $\triangle T_h = T_{\max}-T_{\min}$, which
are the main two parameters controlling instabilities in this problem. 
The system evolves according to the momentum and mass balance equations   
and
to the energy conservation principle. In the equations governing the  
system $u_{r},$ $u_{\phi}$ and $u_{z}$ are the components of the  velocity
field  
 $u$
of the fluid, $T$ the temperature, $p$ the pressure, ${\bf r}$ the  radio vector and $t$ the time are denoted. The magnitudes are expressed
in dimensionless form  
 after rescaling in
the following form: ${\bf r}^{\prime } ={\bf r}/d,$ $t^{\prime  
} =\kappa t/d^{2},$ $%
u^{\prime } =du/\kappa ,$ $p^{\prime } =d^{2}p/\left( \rho  _{0}\kappa  
\nu
\right) ,$ $\Theta  =\left( T-T_{0}\right) /\triangle T$. Here $\kappa$  is the
thermal diffusivity, $\nu $ the kinematic
viscosity of the liquid and $\rho_{0}$ is the mean 
density at the environment temperature  $T_{0}$.\par
The governing dimensionless equations (the
primes in the corresponding fields have been dropped) are the continuity
equation, 
\begin{equation}
\nabla \cdot u =0. \label{1general}
\end{equation}
 The  energy balance equation,
\begin{equation}
\partial _{t}\Theta +u\cdot \nabla \Theta  =\nabla^2 \Theta,  
\label{3general}
\end{equation}
The Navier-Stokes equations, 
\begin{equation}
\partial _{t}u+\left( u\cdot \nabla \right) u =Pr\left( - \nabla  p+\nabla^2  
u+\frac{R\rho }{\alpha \rho _{0}\triangle T}e_{z}\right) , \label{2general}
\end{equation}
where the operators and fields are expressed in cylindrical coordinates  \cite{poz}
and $e_{z}$ is the unit vector in the $z$ direction. 
Here the Oberbeck-Bousinesq  approximation has been used. It consist of  considering
  only in the buoyant term, the
following density dependence on
temperature  $\rho  =\rho_{0}\left[ 1-\alpha \left( T-T_{0}\right)  \right]$,
where $\alpha $ is the thermal
expansion coefficient. The following dimensionless numbers have been  introduced: 
\begin{equation}
Pr =\frac{\nu }{\kappa },\;\;\;R =\frac{g\alpha \triangle Td^{3}}{\kappa  
\nu },
\end{equation}
where $g$ is the gravity constant. $Pr$ is the Prandtl number which  is  assumed to have a large  value 
and $R$ the Rayleigh number, representative of the buoyancy effect.  

\subsection{Boundary conditions}
We discuss now the boundary conditions (bc).  The top  
surface is flat, which implies the following condition on the  velocity,
\begin{eqnarray}
u_{z} & =&0,\;\;\mbox{on}\;\;z =1. 
\end{eqnarray}
 The variation of the surface tension with temperature is considered: $\sigma  \left( T\right)
 =\sigma _{0}-\gamma \left( T-T_{0}\right) ,$ where $\sigma _{0}$ is  
the surface tension at temperature $T_{0 },$ $\gamma $ is the constant  rate
of change of surface tension with temperature ($\gamma $ is positive for
most current liquids). This effect supplies the Marangoni conditions for  the velocity fields
which in  dimensionless form are,
\begin{eqnarray}
 &&\partial _{z}u_{r}+M\partial _{r}\Theta   = 0,\;\;
\partial_{z}u_{\phi}+\frac{M}{r}\partial _{\phi}\Theta  =0,\;\;\mbox{on}\;\;z =1.
\end{eqnarray}
Here $M =\gamma \triangle Td/\left( \kappa \nu \rho _{0}\right) $ is the
Marangoni number. In our particular problem the  
Rayleigh and the Marangoni numbers are related in the same way as in 
the experiments in Ref. \cite{angel} $M =9.2 \cdot 10^{-8} R/d^2$, so buoyancy 
effects are dominant. 
The remaining boundary conditions correspond to  rigid
walls and are expressed as follows,
 \begin{eqnarray}
u_{r} & =&u_{\phi} =u_{z} =0,\;\;\mbox{on}\;\;z =0, \\
u_{r}  & =&u_{\phi} =u_{z} =0,\;\;\mbox{on}\;\;r =a^*,\;\;r  =a^*+\delta^*,
\end{eqnarray}
where $a^* =a/d$ and $\delta^* =\delta/d$.\par
For temperature 
we consider the dimensionless form of Newton's
law for heat
 exchange at the surface, 
\begin{equation}
\partial _{z}\Theta  =-B\Theta ,\;\;\mbox{on}\;\;z =1,
\end{equation}
where $B$ is the Biot number. At the bottom a  linear
profile is imposed, 
\begin{eqnarray}
\Theta & =& \left(-\frac{r}{\delta^*} + \frac{a}{\delta}\right) \frac{\triangle T_h}{\triangle T} + 1,
\;\;\mbox{on}\;\;z =0,
\end{eqnarray}
while in the lateral walls conducting boundary conditions 
are considered
\begin{eqnarray}
\Theta & =&1,\;\;\mbox{on}\;\;r =a^*, \\
\Theta & =&\left(-1+\frac{a}{\delta}\right)\frac{\triangle T_h}{\triangle T} + 1,\;\;\mbox{on}\;\;r =a^*+\delta^*.
\end{eqnarray}
>From this boundary conditions it is clear how not only
$\triangle T$ defining the Rayleigh number is involved in the problem, but also 
$\triangle T_h$.
Due to the fact that pressure is kept in the equations, 
additional
boundary conditions are needed. They  are obtained by the
continuity equation at $z =1$ and  the  
normal component of the momentum equations 
on $r =a^*,$ $r =a^*+\delta^*$ and $z =0$,  \cite{numerico}.
\begin{eqnarray}
&& \nabla \cdot u =0, \;\;\mbox{on}\;\;z =1,\label{bcespecial1} \\
&& Pr^{-1}\left( \frac{\partial u_{r}}{\partial t}+u_{r}\frac{\partial
u_{r}}{\partial r}+\frac{u_{\phi }}{r}\frac{\partial u_{r}}{\partial
\phi }+u_{z}\frac{\partial u_{r}}{\partial z}-\frac{u_{\phi }^{2}}{r}%
\right) \stackrel{}{=}  \nonumber \\
&&\stackrel{}{=}-\frac{\partial p}{\partial r}+\Delta u_{r}-\frac{u_{r}}{%
r^{2}}-\frac{2}{r^{2}}\frac{\partial u_{\phi }}{\partial \phi },
\;\;\mbox{on}\;\;r =a^*,\;r=a^*+\delta^*,\label{bcespecial1} \\
&& Pr^{-1}\left( \frac{\partial u_{z}}{\partial t}+u_{r}\frac{\partial
u_{z}}{\partial r}+\frac{u_{\phi }}{r}\frac{\partial u_{z}}{\partial
\phi }+u_{z}\frac{\partial u_{z}}{\partial z}\right) \stackrel{}{=} 
\nonumber \\
&&\stackrel{}{=}-\frac{\partial p}{\partial z}+\Delta u_{z}-b+RT,
\;\;\mbox{on}\;\;z=0,\label{bcespecial1}
\end{eqnarray}
where $b=d^{3}g/(\kappa \nu)$, and $\Delta= r^{-1} \partial/\partial r 
(r \partial /\partial r) + r^{-2} \partial^2/\partial \phi^2  + \partial^2 /\partial z^2$.
\section{Basic state}
As soon as the lateral walls are at different temperatures and  
a horizontal temperature gradient is set  at the bottom, a  stationary
convective motion appears in the fluid. In contrast to the classical 
B\'enard-Marangoni problem with uniform heating, here the basic state   
is not conductive, but convective.
In order to calculate it,  for computational convenience, the following  change has been performed: 
$r^{\prime } =2 r /\delta^* -2 a^*/\delta^*-1 $
and $z^{\prime } =2z-1,$ which transforms the domain $\Omega_1  =\left[ a^*,
a^*+\delta^* \right] \times \left[ 0,2\pi \right] \times  \left[0,1\right]$ into 
$\Omega =\left[ -1,1\right] \times \left[ 0,2\pi
\right] \times \left[ -1,1\right] $. After these changes,  and
since the basic state has radial symmetry (there is no  
dependence on $\phi $) the steady state 
equations in the infinite Prandtl number approach 
become (the primes have been dropped), 
\begin{eqnarray}
A\partial _{r}p+G^{2}u_{r} & =&\Delta^{*}u_{r},  \label{1} \\
2\partial _{z}p+b-R\Theta & =&\Delta ^{*}u_{z},  \label{3} \\
Gu_{r}+A\partial _{r}u_{r}+2\partial
_{z}u_{z} & =&0,  \label{4} \\
u_{r}A\partial _{r}\Theta +2u_{z}\partial _{z}\Theta & =&\Delta  ^{*}\Theta ,  \label{5}
\end{eqnarray}
where $\Delta ^{*} =A^{2}\partial _{r}^{2}+GA\partial  
_{r}+4\partial _{z}^{2}$, $b =d^3 g/\kappa \nu$, $A =2  
d/\delta$ and 
$G\left( r\right)  =2d/(2a+\delta +r\delta)$.\par
The boundary conditions are now as follows, 
\begin{eqnarray}
&& u_{z}  =2\partial _{z}u_{r}+MA\partial _{r}T = 2  \partial_{z}\Theta +B 
\Theta  =0,\;\;\mbox{on}\;\;z =1, \label{1bc} \\
&& Gu_{r}+A\partial _{r}u_{r}+2\partial_{z}u_{z} =0,\;\;\mbox{on}\;\;z  =1, \label{2bc}\\
&& u_{r}  =u_{z} =0,\Theta  =(1-r)/2,\;\;\mbox{on}\;\;z =-1,  
\label{3bc}\\
&& 2\partial _{z}p+b-R\Theta  =\Delta^{*}u_{z},\;\;\mbox{on}\;\;z  =-1, \label{4bc}\\
&& u_{r}  =u_{z} =0,\,\Theta  =0, \;\;\mbox{on}\;\;r =1,
\label{5bc} \\
&&A\partial _{r}p+G^{2}u_{r} =\Delta^{*}u_{r},\;\;\mbox{on}\;\;r =1,
\label{5bcbis} \\
&& u_{r}  = u_{z} =0,\,\Theta  =1,\;\;\mbox{on}\;\;r =-1,  \label{6bc}\\
&& A\partial_{r}p+G^{2}u_{r} =\Delta ^{*}u_{r},\;\;\mbox{on}\;\;r  =-1. \label{6bcbis}
\end{eqnarray}
\subsection{Numerical method}
We have solved numerically Eqs. (\ref{1})-(\ref{5}) together with the  boundary conditions
(\ref{1bc})-(\ref{6bcbis}) by a Chebyshev collocation method. This  approximation is given by four
perturbation fields $u_{r}\left( r,z\right) ,$ $u_{z}\left(  r,z\right),$ 
$p\left( r,z\right) $ and $\Theta \left( r,z\right) $ which are expanded   
in a
truncated series of orthonormal Chebyshev polynomials:
\begin{eqnarray}
u_{r}\left( r,z\right)  
& =&\sum_{n =0}^{N}\sum_{m =0}^{M}a_{nm}T_{n}\left(r\right)  T_{m}\left( z\right),  \label{a1} \\
u_{z}\left( r,z\right)  
& =&\sum_{n =0}^{N}\sum_{m =0}^{M}b_{nm}T_{n}\left(r\right)  T_{m}\left( z\right),  \label{a2} \\
p\left( r,z\right)  
& =&\sum_{n =0}^{N}\sum_{m =0}^{M}c_{nm}T_{n}\left( r\right)
T_{m}\left( z\right),  \label{a3} \\
\Theta \left( r,z\right)  
& =&\sum_{n =0}^{N}\sum_{m =0}^{M}d_{nm}T_{n}\left( r \right)  T_{m}\left( z\right).  \label{a4}
\end{eqnarray}
Expressions (\ref{a1})-(\ref{a4}) are replaced into the
equations (\ref{1})-(\ref{5}) and boundary conditions 
(\ref{1bc})-(\ref{6bc}). The $N+1$ Gauss-Lobato points 
$\left( r_{j} =\cos (\pi (1 - j/N)),\; j =0,...,N \right) $ in the  $r$ axis  
and the 
$M+1$ Gauss-Lobato points $\left( z_{j} =\cos( \pi (1 - j /M),\;
j =0,...,M \right) $ in the $z$ axis are calculated. The previous
equations are evaluated at these points according to the rules explained   
in Ref. \cite{cedya}, in this way $4\left( N+1 \right) \left( M+1 \right)$  equations
are obtained with  
 $4\left(N+1\right) \left( M+1\right) $ unknowns. The system has not  
 maximun
rank because pressure is only determined up to
a constant value. Since this value is not affecting to the other  physical
magnitudes, we replace the evaluation of the normal component of  
 the momentum equations at ($r_{j=N}=1,$ $z_{j=4} =\cos  (\pi(1-4/M))$) by 
a value for the pressure at  
 this point, for instance $p=0$.
 To solve the resulting nonlinear equatins,  a Newton-like  iterative method is used.
In a first step the nonlinearity is discounted and a solution is found  by solving
the linear system. It is corrected by small perturbation  fields: $\bar{u}_{r}\left( r,z\right) ,$
$\bar{u}%
 _{z}\left( r,z\right) ,$ $\bar{p}\left( r,z\right) $ and  
$\bar{\Theta}\left(
r,z\right) ,$
\begin{eqnarray}
u_{r}^{i+1}\left( r,z\right) & =&u_{r}^{i}\left( r,z\right) +\bar{u}%
_{r}\left( r,z\right) , \\
u_{z}^{i+1}\left( r,z\right) & =&u_{z}^{i}\left( r,z\right) + \bar{u}
_{z}\left( r,z\right) , \\
p^{i+1}\left( r,z\right) & =&p^{i}\left( r,z\right) +\bar{p}\left(  
r,z\right)
, \\
\Theta ^{i+1}\left( r,z\right) & =&\Theta ^{i}\left( r,z\right)  
+\bar{\Theta}%
\left( r,z\right),
\end{eqnarray}
in such way that solutions at $i+1$ step are obtained after solving
Eqs. (\ref{1})-(\ref{6bc}) 
linearized around the approach at step $i$. The considered criterion of
convergence is that the difference between two consecutive  
approximations in $l^{2}$ norm should be  smaller than $10^{-9}.$
\subsection{Results on the basic state}
Previous theoretical works \cite{mercier,parmentier,smith} 
dealing with this problem have approached the basic state solution
as a  parallel flow, where only the horizontal component 
of the velocity field exists. This approach requires the imposition of
a constant temperature gradient over all the fluid layer and
at the top boundary. In those descriptions
this gradient becomes the main control parameter.  
In our model lateral effects are considered,
so parallel flow is no longer valid and no restrictions are
needed in the temperature boundary conditions which, as explained
in the previous section, is the usual Newton law. This boundary condition
is the origin of two control parameters related to temperature, i. e.,
$\triangle T_h$ and $\triangle T$. This possibility, which had been already addressed
in \cite{pof}, is explored in detail in the next section.

We have found two types of basic states. The first one, displayed 
in Fig. 2 shows the isotherms obtained at $B=1.25$, $\delta^*=10$ and
 $\Delta_h T=\triangle T=6.4^{\rm o} \,{\rm C}$. 
The bottom profile is linear, as expected from the boundary condition, however
at the top it  approaches the constant ambient temperature, since the Biot number is large. 
This field is similar to that of 
the layer heated from below, i.e. roughly
speaking hotter at the bottom. This feature is shared with the {\it linear flow}
described in \cite{smith}. Fig. 3 shows the velocity fields for 
the same conditions. They are formed by two co-rotative 
rolls perpendicular to the gradient. This result is  similar to that
obtained in experiments reported in Refs. \cite{angel,daviaud,burguete}, where 
co-rotative rolls perpendicular to the gradient are found as well. 
These  results coincide with  those obtained in Ref. \cite{pof} 
for a rectangular geometry. \par

The second kind of basic state is represented in
Fig. 4 for $B=0.8$, $\delta^* = 2.5$, 
$\Delta_h T = 0.3^{\rm o} \,{\rm C}$, $\triangle T = 1.84^{\rm o} \,{\rm C}$. 
The most striking feature is that along the vertical axis there are changes in the sign of the temperature gradient. It seems that  the Biot number
is not  dissipating heat effectively at the top and then heat 
is advected by the velocity field. This temperature profile coincides with that of the {\it return flow}  defined in Ref. \cite{smith}, where 
the role of the transport of heat by the  velocity field as the origin of the instability is discussed.
As we will demonstrate in the next section  this flow becomes unstable through 
 an oscillatory bifurcation, as it is also  reported in \cite{smith}. 
 Crucial to the appearence of the two different flows described in this section, is not only the Biot number that dissipates heat at the top,
but also the role of parameters $\triangle T$ and $\triangle T_h$ 
which introduces heat into the system. 

Fig. 5 shows a transition between the two states we have just described.
It is obtained at $B=0.3$, $\delta^* = 10$, 
$\Delta_h T = 10^{\rm o} \,{\rm C}$, $\triangle T = 20.41^{\rm o} \,{\rm C}$. At this stage the flow already becomes unstable through 
 an oscillatory bifurcation.

\section{Linear stability of the basic state}
The stability of the basic state is studied by perturbing it
 with a vector field depending on the 
$r,\phi $ and $z$ coordinates, in a fully 3D analysis: 
\begin{eqnarray}
&&u_{r}\left( r,\phi ,z\right)  =u_{r}^{b}\left( r,z\right)  +\bar{u}%
_{r}\left( r,z\right) e^{i{\rm m}\phi+ \lambda t },  \label{10} \\
&&u_{\phi }\left( r,\phi ,z\right)  =u_{\phi }^{b}\left(  
r,z\right) +
\bar{u}_{\phi }\left( r,z\right) e^{i{\rm m}\phi +\lambda t},  \label{20}  \\
&&u_{z}\left( r,\phi ,z\right)  =u_{z}^{b}\left( r,z\right)  +\bar{u}%
_{z}\left( r,z\right) e^{i{\rm m}\phi+\lambda t },  \label{30} \\
&&\Theta \left( r,\phi ,z\right)  =\Theta ^{b}\left( r,z\right)  
+\bar{\Theta
}\left( r,z\right) e^{i{\rm m}\phi+\lambda t },  \label{40} \\
&&p\left( r,\phi ,z\right)  =p^{b}\left( r,z\right) +\bar{p}\left(
r,z\right) e^{i{\rm m}\phi +\lambda t}.  \label{50}
\end{eqnarray}
Here the superscript $b$ indicates the
corresponding quantity in the basic state and the bar refers to the 
perturbation. We have considered Fourier modes expansions
 in the angular direction, because along it  the boundary conditions 
are periodic. 
We replace the expressions (\ref{10})-(\ref{50}) 
into the basic equations and after linearizing the resulting system, we
obtain the following eigenvalue problem (the bars have been  
dropped): 
\begin{eqnarray}
\Delta_{{\rm m}}u_{r}-A\partial _{r}p-G^{2}u_{r}-2G^{2} i {\rm m} u_{\phi }  
& =&0,  \label{eig1} \\
\Delta _{{\rm m}}u_{\phi }-Gi{\rm m}p+2G^{2}i{\rm m}u_{r}-G^{2}u_{\phi } & =&0,   
\label{eig2}
\\
\Delta_{{\rm m}}u_{z}-2\partial _{z}p+R \Theta & =&0,  \label{120} \\
Gu_{r}+A\partial _{r}u_{r}+Gi{\rm m}u_{\phi }+2\partial _{z}u_{z} & =&0,
\label{eig3} \\
\Delta _{{\rm m}} \Theta-u_{r}A\partial _{r} \Theta^{b}-u_{r}^{b}A\partial  
\Theta-2u_{z}^{b}\partial
_{z} \Theta -2u_{z}\partial _{z} \Theta^{b} & =&\lambda \Theta,   
\label{eig4}
\end{eqnarray}
where $\Delta_{\rm m}  = A^2 \partial_r^2 + GA \partial_r - {\rm m}^2 G^2 + 4  
\partial_z^2$. The following boundary conditions for the pertubations  are obtained,
\begin{eqnarray}
u_{z} & =&2\partial _{z}u_{r}+MA\partial _{r}\Theta  =2\partial  
_{z}u_{\phi
}+GM i {\rm m} \Theta  =2 \partial _{z}\Theta +B\Theta  =0, 
\;\;\mbox{on}\;\;z =1,  \label{bceig1} \\
u_{r} & =&u_{\phi } =u_{z} =\Theta  =0,\;\;\mbox{on}\;\;z  =-1,   
\label{bceig2} \\
u_{r} & =&u_{\phi } =u_{z} =0,\; \Theta  =0,\;\;\mbox{on}\;\;r  =1,   
\label{bceig3} \\
u_{r} & =&u_{\phi } =u_{z} =,\; \Theta  =0,\;\;\mbox{on}\;\;r  =-1,   
\label{bceig4}
\end{eqnarray}
together with 
\begin{eqnarray}
&& \Delta_{\rm m} u_r - A\partial _{r}p-2G^{2}i {\rm m} u_{\phi  
} =0,\;\;\mbox{on}\;\;r =\pm 1,  \label{bceig34bis} \\
&& \Delta_{\rm m} u_z-2\partial _{z}p+R\Theta  =0,\;\;\mbox{on}\;\;z =-1,    
\label{bceig2bis} \\
&&Gu_{r}+A\partial _{r}u_{r}+Gi{\rm m}u_{\phi }+2\partial  
_{z}u_{z} =0,\;\;\mbox{on}\;\;
z =1.  \label{eig1bis}
\end{eqnarray}
\subsection{Numerical method}
The  eigenvalue problem is discretized 
with the  Chebyshev collocation method used for the basic state. 
Now there is a new field, the angular velocity, which is expanded as  follows, 
\begin{equation}
u_{\phi }(r,z) =\sum_{n =0}^{N}\sum_{m =0}^{M}e_{nm}T_{n}\left(   
r\right)
T_{m}\left( z\right).
\end{equation}
In order to calculate the eigenfunctions and thresholds 
of the generalized eigenvalue problem the equations are posed at the  collocation points 
according to the rules explained in \cite{cedya}, so that  a total 
of $5 (N+1) (M+1)$ algebraic 
equations are obtained with the same number of unknowns. If the 
coefficients of the unknowns which 
form the matrices $A$ and $B$ satisfy 
det$(A- \lambda B) =0$, 
a nontrivial solution of the linear homogeneous 
system exists. This condition generates a dispersion 
relation $\lambda \equiv \lambda({\rm m},R,M,B,u^b,\Theta^b,p^b)$, 
equivalent to calculate directly the eigenvalues from the 
system $AX =\lambda BX$, 
where $X$ is the vector which contains the 
unknowns. If Re$(\lambda)<0$  the basic state is 
stable while if Re$(\lambda)>0$ the basic state becomes unstable. 
In the Appendix a study of the convergence of
the thresholds calculated with this method is shown.
  
\subsection{Results on the stability analysis}
The basic solutions obtained in section 3, become unstable when the control parameter, $\triangle T$ in this work, is increased. When  $\triangle T$ is changed, $\triangle T_h$ which is not the control parameter in our choice, is fixed to a non zero value. In  the explored parameter range two types of bifurcations take place: 
stationary and oscillatory. Figure 6 depicts the m dependence  of 
Re$(\lambda)$ for $B=1.25$, $\delta^* =10$ and $\triangle T_h = \triangle T$ at the 
threshold ($\triangle T_c =6.4^{\rm o} \,{\rm C}$) and below it 
($\triangle T < 6.4^{\rm o} \,{\rm C}$). 
The bifurcation is
stationary since  the imaginary part of the eigenvalue $\lambda$ is
 zero  at the critical value m$_c =20$. 
Figure 7 is similar to Figure 6 but at 
$B=0.5$, $\delta^* = 10$ and $\triangle T_h = 5^{\rm o} \,{\rm C}$ and 
$\triangle T_c = 8.63^{\rm o} \,{\rm C}$. Now the bifurcation is oscillatory since  
the imaginary part of the eigenvalue is non zero at the critical value m$_c =18$. 
In order to understand how physical conditions affect the instability,  
Fig. 8 summarizes the influence of heat related
parameters in the bifurcation. There, it is shown the critical 
value of $\triangle T_c$, as a function of the Biot number $B$, for several 
values of $\triangle T_h$. In this figure stationary and oscillatory 
bifurcations are observed and it shows that having one or  other 
 transition depends on an equilibrium between these
parameters. Travelling waves are favoured by systems storing a lot of heat, 
i.e., low values of the Biot number, and high values of $\triangle_h T$, which 
facilitates the transport of heat by the velocity field \cite{smith}. 
We notice here the novelty of these results, 
where transitions between stationary rolls and
oscillatory waves are due to changes in heat related parameters. 
This has not been reported in any experimental or theoretical work before, where transitions have been found due to variations in Prandtl number \cite{smith}
or in aspect ratio \cite{daviaud,burguete}.

In Fig. 8 it is also clear that thresholds decrease as $B$ increases 
while the dependence of thresholds on
$\triangle T_h$ is not monotonous, it presents a minimum for $\triangle T_h = 2^{\rm o} \,{\rm C}$ 
as figure 9 a) displays. 
In Fig. 9 b) it is shown how the critical wavenumber  
increases with $\triangle T_h$, on the other hand it remains almost 
constant with $B$. \par
The geometry and size of the box is also affecting
the instability as it is noticed mainly in low aspect ratio systems.
Although a detailed  study of patterns in small containers
is beyond the scope of this work,  table I gives an insight into it.
It shows the critical 
value of $\triangle T_c$ and the corresponding critical wavenumber m$_c$ 
for some $\delta^*$. The wavenumber does not change monotonously
with the aspect ratio, as also happens in small containers
with uniform heating from below \cite{dauby2,numerico,mluisa} as well. \par
The critical wavenumber for the stationary bifurcation 
shows that a 3D structure appears which 
consists of radial rolls whose axes are 
parallel to the radial coordinate as 
Fig. 10 shows for $B=1.25$ and $\delta^*=10$ and $\triangle T_c=6.4^o$C. 
The growing perturbation along the transverse plane
 is plotted in Fig. 11, where a structure appears near
 the hot side as in the experiment reported in  \cite{de}.
This is so because, in this case, the instability has the same origin
as a static fluid heated from below, so since the  
 vertical temperature gradient is greater at this boundary, perturbation
grows there first. On the other hand since the hot cylinder  counteracts the vertical temperature gradient it inhibits 
the instability threshold. 
These results are in good agreement with those obtained in Ref. \cite{pof} in
similar conditions but in a cartesian geometry.  In that article 
at a threshold  $\triangle T_c =6.8^o$C, the critical wavenumber 
is $k_c=2.4$, which would correspond to  $12-36$ periods respectively in
the inner and outer radii of an annulus as ours. 
This matches rather well with our critical threshold and  wavenumber 
which is ${\rm m}_c =20$, suggesting that, at least
for this aspect ratio, geometry does not strongly affect the results.
These results and discussions of the stationary bifurcation
are in  good agreement with those reported in Refs.  
\cite{daviaud,mercier,parmentier,smith,gershuni}. \par
We also obtain oscillatory bifurcations.  For instance, there is one 
at $\delta^*=10$ ($d=2$ mm), $\triangle T_h = 5^o$C, $B=0.5$ and $\triangle T_c=8.63^o$C, 
which corresponds to travelling waves or  hydrothermal waves described 
in Refs. \cite{daviaud,angel,burguete}. Fig. 12 depicts 
the eigenfunction along the transverse plane, where 
the tilted axis of the waves can be seen. This result is comparable
to that reported in Ref. \cite{burguete} where 
for a cartesian container with $L_x =20$ mm, $d \sim 2$ mm, hydrothermal waves of 
wavenumber $k_c \sim 1$, emerge above a threshold.
This would correspond to a pattern with
 $6-19$ periods in the inner and outer radii of our annulus respectively,
which is not far from our critical wavenumber m$_c=18$. 
The angle of inclination of the waves in Fig. 10 is $ \phi \sim 3 \pi/4$ 
also  close to that of  $2.60-1.75$ rad of  Ref. \cite{burguete}. 

\section{Conclusions}
We have studied a BM lateral heating problem in a cylindrical annulus.  
The problem has been solved with a Chebyshev   
collocation method  by keeping the original Navier-Stokes
equations, where appropriate   boundary conditions are required for pressure \cite{numerico}. The scheme developed in Ref. \cite{numerico} for stability problems has been extended  to calculate a non trivial basic state. 
The procedure is confirmed to be reliable,  effective and easy to implement.\par

We have obtained two kinds of basic solutions.
One type is formed by co-rotative rolls such as those reported in experiments \cite{angel,daviaud,burguete}, and in similar theoretical studies \cite{pof}. The second kind of solution has similar features to those of 
the {\it return flow} described in  \cite{smith}.

The linear stability analysis of these solutions shows stationary 
bifurcations to radial rolls and oscillatory bifurcations 
to hydrothermal waves, whose appearance has been proved to  depend on
their heat related parameters ($B$, $\triangle T_h$ and $\triangle T$). This
fact has not been adressed before, perhaps because experimentally these
parameters are hard to manipulate.  Properties of the growing patterns
 coincide with those reported in previous works.
Stationary solutions appear in basic flows with temperature profiles close to
 the static fluid heated from below \cite{pof, smith}. Their 
wavenumbers are in the range of experimental
 \cite{angel,daviaud,burguete} and theoretical studies \cite{pof}.
Also the structure in the $z-r$ plane agrees with 
previous results \cite{pof,de}.
Oscillatory solutions emerge from a basic state of the type of {\it  return flow}
as predicted in \cite{smith}. Their structure is comparable to that
obtained in experiments \cite{burguete}.  \par
\medskip
\noindent
{\bf Acknowledgments}\par
\medskip
We gratefully thank Christine Cantell, Maureen Mullins and 
Angel Garcimart\'{\i}n for useful comments and suggestions. 
This work was partially supported by a Research Grant  
MCYT (Spanish Government) BFM2000-0521 and by the University 
of Castilla-La Mancha.

\newpage
\section*{Appendix}
\subsection*{ Convergence of the global numerical method} 
To carry out a test on the convergence of the global numerical method 
(calculation of the basic state and subsequent linear stability  
analysis) we
compare the differences in the thresholds of the differences of  
temperature
($\triangle T$) to different orders of expansions for some values of 
the parameters involved in the problem: 
$\delta^* \in [8.26,10.52]$. 
In table II the thresholds for these states 
varying the aspect ratio are shown for four consecutive
expansions varying the number of polynomials taken in the 
$r$ ($N$) and  $z$ ($M$) coordinates. 
>From the table it is seen that convergence is reached 
within a relative precision for $\triangle T_c$ of $10^{-2}$. 
If $M$ is increased  the difference between successive 
expansions is order $10^{-1}$ while increasing $N$ it is 
order $10^{-2}$. Further increments of the order of the expansions give
thresholds values oscillating around the quoted ones, 
so within that error, thresholds are convergent. 
The range of parameters has been selected in order to reach this 
convergence criterium. 
\newpage
\noindent
{\bf Table captions}\par
\bigskip
\noindent
{\bf Table I}\par
\smallskip
\noindent
Critical $\triangle T$ and critical m$_c$ vs the aspect ratio $\delta^*$.  
The value of the horizontal temperature difference is $\triangle T_h=0.4$ 
and $B=1$.\par
\bigskip
\noindent
{\bf Table II}\par
\smallskip
\noindent
Critical temperature differences ($^{\circ}{\rm C}$)  for different  
values 
of the aspect ratio at consecutive orders in the expansion in Chebyshev   
polynomials
($B =1.25$).\par
\bigskip
\noindent
{\bf Figure captions}\par
\bigskip
\noindent
{\bf Figure 1}\par
\smallskip
\noindent
Problem set up ($a=0.01$ m, $\delta=0.02$ m).\par
\bigskip
\noindent
{\bf Figure 2}\par
\smallskip
\noindent
Isotherms of the basic state 
corresponding to values of  the parameters 
$R =2228,\, M =51,\, \delta^* =10,\;B =1.25$ ($\Delta T= 6.4^o$C, 
$\Delta T_h = 6.4^o$C, $d=2$ mm). \par
\bigskip
\noindent
{\bf Figure 3}\par
\smallskip
\noindent
 Velocity field at same conditions of Fig. 2.\par
\bigskip
\noindent
{\bf Figure 4}\par
\smallskip
\noindent
Isotherms and velocity field of the basic state 
corresponding to  values of the parameters 
$R = 40995,\, M = 59,\, \delta^* = 2.5,\;B = 0.8$ ($\Delta T= 1.84^o$C, 
$\Delta T_h = 0.3^o$C, $d=8$ mm). \par
\bigskip
\noindent
{\bf Figure 5}\par
\smallskip
\noindent
Isotherms and velocity field of the basic state 
corresponding to values of the  parameters 
$R = 7105,\, M = 163,\, \delta^* = 10,\;B = 0.3$ ($\Delta T= 20.41^o$C, 
$\Delta T_h = 10^o$C, $d=2$ mm). \par
\bigskip
\noindent
{\bf Figure 6}\par
\smallskip
\noindent
Maximum real part of the growth rate $\lambda$ as a function of ${\rm m}$ for   
basic state at $M =51$, $\delta^* =10$, $B =1.25$. 
The top line stands for the threshold condition $R_c =2228$ (or $\triangle T_c =6.4^o$C), 
and its maximum determines the critical 
${\rm m}_c$. The bottom line is for 
$R <  R_c$. Solid lines correspond to branches with 
real eigenvalues while the dotted ones  stand for complex eigenvalues.\par
\bigskip
\noindent
{\bf Figure 7}\par
\smallskip
\noindent
Maximum real part of the growth rate $\lambda$ as a 
function of ${\rm m}$ for   
a basic state at $M = 69$, $\delta^* = 10$, $B =0.5$, at
the threshold condition $R_c = 3004$ (or $\triangle T_c = 8.63^o$C).
The maximum determines the critical ${\rm m}_c$. 
The solid line correspond to a branch with real eigenvalues while the dotted one  
is for complex eigenvalues.\par
\bigskip
\noindent
{\bf Figure 8}\par
\smallskip
\noindent
Critical $\triangle T$ vs the Biot number $B$ for different values of 
$\triangle T_h$ ($\delta^*=10$).\par
\bigskip
\noindent
{\bf Figure 9}\par
\smallskip
\noindent
Critical $\triangle T$ and m$_c$ vs $\triangle T_h$ for $\delta^*=10$ 
and $B=0.7$.\par
\bigskip
\noindent
{\bf Figure 10}\par
\smallskip
\noindent
Isotherms of the growing perturbation in the $x-y$ plane at a stationary   
bifurcation point 
  $(\triangle T = \triangle T_h =6.4^{\circ}\, {\rm C},\, \delta^* =10,\, B =1.25)$. \par
\bigskip
\noindent
{\bf Figure 11}\par
\smallskip
\noindent
a) Isotherms  of the growing perturbation in the $r-z$ plane at a  
stationary bifurcation point  
  $(\triangle T = \triangle T_h =6.4^{\circ}\, {\rm C},\, \delta^* =10,\, B =1.25)$. 
b)  Velocity field of the growing perturbation in the $r-z$ plane at  
same conditions.\par
\bigskip
\noindent
{\bf Figure 12}\par
\smallskip
\noindent
a) Isotherms of the eigenfunction in the $x-y$ plane at an 
oscillatory bifurcation point  
$(\triangle T_h = 5^{\circ}\, {\rm C},\,\triangle T = 8.63^{\circ}\, {\rm C},\, \delta^* =10,\, B =0.5)$. \par
\bigskip
\noindent
{\bf Figure 13}\par
\smallskip
\noindent
a) Isotherms of the eigenfunction in the $r-z$ plane at the  
oscillatory bifurcation point considered in Fig. 12. 
b)  Velocity field of the eigenfunction in the $r-z$ plane at  
same conditions.\par
\bigskip
\noindent
\newpage
\centerline{\bf Table I}
\begin{center}
\begin{tabular}{|c|c|c|} \hline 
$\delta^* $ & $\triangle T_c $ & ${\rm m}_c $ \\ 
\hline
$2.5$ & $3.91$ & $11$  \\ 
$10/3$ & $2.63$ & $1$  \\ 
$5$ & $3.28$ & $2$ \\ 
$10$ & $4.03$ & $10$ \\ 
\hline
\end{tabular}
\end{center}
\centerline{\bf Table II}
\begin{center}
\begin{tabular}{|c|c|c|c|c|} \hline 
$\delta^* $ & $13\times 9$ & $13\times 13$ & $25\times 13$ & $29\times 13$   \\ 
\hline
$10.52$ & $6.8711$ & $6.8519$ & $6.9551$ & $6.9229$ \\ 
$10.25$ & $6.5131$ & $6.5224$ & $6.6792$ & $6.6665$ \\ 
$10.00$ & $6.2223$ & $6.2298$ & $6.3717$ & $6.3938$ \\ 
$9.52$ & $5.8794$ & $5.9490$ & $5.9126$ & $5.9519$ \\ 
$9.09$ & $5.8415$ & $6.0807$ & $5.6626$ & $5.6868$ \\ 
$8.26$ & $5.7341$ & $5.9204$ & $5.6812$ & $5.6892$ \\
\hline
\end{tabular}
\end{center}

\end{document}